\documentclass{article}

\usepackage{ifxetex}
\usepackage{commath}
\usepackage{amsthm}
\usepackage{amssymb}
\usepackage{enumitem}

\newtheorem{mydef}{Definition}
\ifxetex
  \usepackage{fontspec}
\else
  \usepackage[T1]{fontenc}
  \usepackage[utf8]{inputenc}
  \usepackage{lmodern}
\fi
\newtheorem{theorem}{Theorem}
\title{The smallest primitive root modulo a prime}
\author{Jana Pretorius\\ School of Physical, Environmental and Mathematical Sciences \\ The University of New South Wales Canberra, Australia \\ u5783790@anu.edu.au}

\begin{document}
\maketitle

\begin{abstract}
In this paper we will consider new bounds on the smallest primitive root modulo a prime. we will make more judicious use of the P\'olya--Vinogradov and Burgess inequalities, and use them to prove that the smallest primitive root is smaller than $p^{0.68}$ for all primes $p$.
\end{abstract}

\section{Introduction}

This paper will detail new research on upper bounds on the smallest primitive root modulo a prime, which will be referred to at $g(p)$. In \cite{B}, Burgess showed that the smallest possible bound is
\begin{equation}\label{sweep}
g(p) \ll p^{1/4+\epsilon} \quad \textrm{for any} \quad \epsilon>0.
\end{equation}

Since any case of $\epsilon< \frac{1}{4}$ has been too difficult to prove to date, some previous research has focused on the Grosswald conjecture,

\begin{equation*}
g(p)< \sqrt{p}-2,
\end{equation*}

which Cohen, Oliveira e Silva and Truigian proved for all $409<p<2.5\times10^{15}$ and $p>10^{71}$ (Theorem 1.1 in \cite{COT}). Others, including Hunter, have instead attempted to find the lowest value of $\epsilon$ for which (\ref{sweep}) holds for all primes. This paper will  lower the universal bound.

Current literature has focused primarily on using the P\'olya--Vinogradov inequality and Burgess bound on Dirichlet characters separately, but never in conjunction. We will be using both in order to minimise the bounds further. The smallest power $\alpha$ of $p$ which could possibly be proven for all primes is $\alpha =\log2/\log3 = 0.63093...$ as this corresponds to the smallest primitive root of $3$, which is $2$.

This paper will prove two main theorems.
\begin{theorem}
Let $g(p)$ denote the least primitive root modulo $p$, prime. Then \begin{equation}\label{single} g(p)<p^{0.6309}\quad \textrm{for all}\quad p>2.67\times10^{32}. \end{equation}
\end{theorem}
The first theorem takes the power of $p$ slightly smaller that the minimum possible ($\log2/\log3$), $0.6309$,  and identifies the range of $p$ over which we can prove that the bound holds. As we cannot prove that this holds for all $p$, we have found the lowest exponent which can be proven for all $p$, given in Theorem 2.
\begin{theorem}
\begin{equation*}
g(p)<p^{0.68}\quad \textrm{for all primes}\quad p.
\end{equation*}
\end{theorem}

This new bound represents significant improvement on previous findings. Hunter proved that the bound on the least \textit{square-free} primitive root was $p^{0.88}$ for all primes $p$, and showed that $\alpha=0.6309$ holds for all $p> 9.63\times 10^{65}$ (Theorems 0.1 and 0.2 in \cite{H}).

These bounds will be shown using the same methods employed by Hunter. In \S2 and \S3, we will the P\'olya--Vinogradov and Burgess inequalities to establish over which values of $p$ each bound is most useful. We will also identify the ranges of $p$ over which various values of $\alpha$ hold. A sieving inequality will be used in \S5 to tighten these bounds further. Finally, in \S6, computations using a prime divisor tree will provide some small improvements.

\section{Comparing the P\'olya--Vinogradov and Burgess inequalities}

The first step towards lowering the bound on primitive roots is to use improved bounds on the P\'olya--Vinogradov and Burgess Inequalities. Both inequalities place an upper bound on the sum on $\chi (p)$, the non-principal Dirichlet characters modulo $p$, defined as
\begin{equation*}\label{bound}
{S_H(N)} = \abs{\sum_{n=N+1}^{N+H}\chi(n)} \quad \text{for some} \quad H<p.
\end{equation*}
Theorem 2.1 in \cite{COT}, quoting Trevi\~no \cite{T} gives the Burgess bound:
\begin{equation}\label{burgess}
C(r)H^{1-\frac{1}{r}}p^{\frac{r-1}{4r^{2}}}(\log{p})^{\frac{1}{2r}},
\end{equation}
where $r$ is an integer of our own choosing. Below, we will discuss which value of $r$ returns the tightest bound. Hunter gives the P\'olya--Vinogradov inequality as follows (Theorem 1.2 in \cite{H}):
\begin{equation*}\label{PV}
c\sqrt{p}\log{p},
\end{equation*}
where $c$ is a constant. The Burgess bound (\ref{burgess}), unlike P\'olya--Vinogradov, is dependent on the length of the sum, $H$.

In order to compare these two inequalities, it is first necessary to determine what value of $r$ to use in the Burgess bound. The values of $C(r)$ were taken from \S1 of Trevi\~no for the case of $H=10^{15}$ \cite{T}. This value of $H$ was chosen as this was the closest value to $2.5\times10^{15}$, below which Theorem 1.1 in \cite{COT} holds.

By substituting $H=p^\alpha$, we can compute that as $\alpha$ increases, $r=2$ gives the smallest bound for an increasing range of primes. At $\alpha=0.7$, $r=2$ is the preferred value of $r$ for all $p>1.5\times 10^{6}$.

Since $r=2$ provides the smallest Burgess bound within the desired range of $\alpha$'s under consideration, this value will be used to compare the Burgess and P\'olya--Vinogradov bounds, by solving:
\begin{equation*}\label{kimmyb}
C(2)p^{\frac{\alpha}{2}+\frac{3}{16}}(\log p)^{\frac{1}{4}} \leq cp^{\frac{\alpha}{2}}\log p.
\end{equation*}

From this, it follows that for values of $\alpha>\frac{5}{8}$, the P\'olya--Vinogradov inequality (on the right hand side of the above equation) will always give a smaller bound than Burgess. i.e.\ The equation is not true for any $p$. Below $\frac{5}{8}$, Burgess becomes the better bound for a range of $p$ dependent on $\alpha$. For example, at $\alpha=0.6$, Burgess gives a smaller bound for all $p>10^{22}$.

From this, it is clear that for the range of $\alpha$ under consideration in this paper, the P\'olya--Vinogradov bound is the better choice. recall that the smallest $\alpha$ which may be proven for all primes is $\frac{\log2}{\log3}>\frac{5}{8}$.

\section{Results without sieving}

Using the P\'olya--Vinogradov inequality, we now make use of improvements to the constant, $c$. While Hunter used $c=\frac{1}{2\pi}+\frac{1}{\log p}+\frac{1}{\sqrt p \log p}$, Frolenkov and Soundararajan found an improved constant (\S4 of \cite{FS}):

\begin{equation*}\label{redgrave}
c(p)=\frac{1}{2\pi}+\frac{0.8203}{\log p}+\frac{1.0284}{\log p\sqrt p}.
\end{equation*}

Frolenkov and Soundararajan, however, optimised their constant to be the smallest over all primes. For the purposes of this research, it is better to optimise the constant for $p=2.5 \times 10^{15}$. Following the derivation in \S4 of \cite{FS}, we chose a generic value for the constant $L$ of the form $\alpha q^{\beta}$ and find the values of each variable which give the minimum bound. From this, we ascertain that the smallest, and hence best, $c$ to use is:

\begin{equation*}
c(p)=\frac{1}{2\pi}+\frac{1}{\pi \log p}\left(0.4325+\frac{10.15 + \sqrt{p}}{1+ \sqrt{p}}+\frac{1}{\sqrt{p}}\right).
\end{equation*}

Since $c(p)$ is slowly decreasing in $p$, we can set $c=c(p_0)$ which will hold for all $p>p_0$ Following the derivation in \S2 of \cite{H}, we sum the following indicator function over $x$.
\begin{mydef} given a prime $p$, and provided that  $p\nmid x$, we can define the function \begin{equation}\label{pair}f(x)=\frac{\phi(p-1)}{p-1}\left\{1+\sum_{d|p-1,d>1}\frac{\mu(d)}{\phi(d)}\sum_{\chi_d}\chi_d(x)\right\} = \begin{cases}&1\, \textrm{if n is a primitive root}\\&0\, \textrm{otherwise,} \end{cases}
\end{equation}
where $\phi$ is Euler's function, $\mu$ is the M\"obius function, and $\chi_d$ are the characters of order $d$ (Lemma 2.1 in \cite{H}).
\end{mydef}

This sum will be greater than $0$ once we have found a primitive root. We  apply the P\'olya-Vinogradov inequality to (\ref{pair}) and proceed according to the derivations in \cite{COT},\cite{H},\cite{CST}. We conclude that a primitive root exists below $g(p)=p^{\alpha}$ if and only if

\begin{equation*}
p^{\alpha} -(2^{\omega (p-1)}-1)c \sqrt{p}  \log p  > 0.
\end{equation*}

We make use of the following estimate for the value of $\omega(n)$ from Robin (Theorem 16 in \cite{R}) to obtain an equation wholly in terms of $p$ and $\alpha$

\begin{equation*}\label{pinsent}
\omega(n) \leq \frac{\log n}{\log\log n} \left(1 + \frac{1}{\log\log n} + \frac{2.8973}{(\log\log n)^2}\right),
\end{equation*}

where $\omega (n)$ denotes the number of distinct prime divisors of $n$. Following this substitution, we compute the value $p$ above which the inequality holds for a range of $\alpha$. We then calculate the corresponding value of $\omega(p)$. For $\alpha=0.7$, the inequality holds for all $p>5\times 10^{1295}$, and so $\omega(n)=554$. From this we calculate the smallest product of primes (the smallest $\omega(n)$) which exceeds $p=5\times 10^{1295}$, knowing that it must be smaller than 554. Any values up to and including this $\omega(n)$ cannot make the inequality hold, and hence must be checked through other means.

Table 1 identifies the values of $w(n)$ which remain to be checked after initial use of the P\'olya--Vinogradov inequality.
\begin{table}[h]
\begin{center}
\caption{Exceptions of $\omega (p-1)$}
\begin{tabular}{||c c c||}
 \hline
 $\alpha$ & Lower bound on $w(p-1)$  & Upper bound on $w(p-1)$ \\ [0.5ex]
 \hline\hline
 0.8 & 5 & 30 \\
 \hline
 0.75 & 5 & 46 \\
 \hline
 0.7 & 5 & 85 \\
 \hline
 0.65 & 5 & 237 \\
 \hline
 0.6309 & 5 & 437 \\ [1ex]
 \hline
\end{tabular}
\end{center}
\end{table}

\section{Improvements using a sieve}

The sieve used in \cite{COT}, \cite{CST}, \cite{H}, and \cite{MTT}  makes use of the number of small primes dividing $p-1$ in order to eliminate most values of $\omega(n)$ identified as possible exceptions in \S 3. The sieve relies on the use of e-free integers:

\begin{mydef} Let $p$ be a prime and let $e$ be a divisor of $p-1$. Suppose $p\nmid n$ then $n$ is \textbf{e-free} if, for any divisor $d$ of $e$, such that $d>1$, $n\equiv y^{d}$ (mod $p$) is insoluble. \end{mydef}

Following Proposition 3.2 in \S3 of \cite{H}, it can be shown that $n$ is a primitive root if and only is it is $(p-1)$ free.

Consequently, an indicator function has been developed to test whether any particular number $n$ is $p-1$ free (\S3 of \cite{H}). \begin{mydef} \begin{equation*}\label{efree} f(x) = \frac{\phi(e)}{e} \sum_{d|e}\frac{\mu(d)}{\phi(d)} \sum_{\chi \in \Gamma_{d}}{\chi (n)}=\begin{cases}&1\quad \textrm{if n is e-free}\\&0\quad \textrm{otherwise.} \end{cases} \end{equation*} \end{mydef}
Again, we apply the P\'olya--Vinogradov bound the indicator function. Following the derivation in \S3 of \cite{H}, we obtain a new inequality (Theorem 3) to solve. This inequality differs from that in \cite{H} as the square-free conditions on the primitive roots are removed.

\begin{theorem}

Let $\delta>0$ be defined as $$\delta = 1-\sum_{i=1}^{s}p_i^{-1}$$ and let
$\Delta$ be defined as $$\Delta = \frac{s-1}{\delta}+2,$$ where $1 \leq s \leq \omega (p-1)$. Then we are able to find a primitive root $g(p)<p^{\alpha}$ if the following inequality holds: $$\frac{p^{\alpha-0.5}}{\log(p)}>c(2^{n-s}\Delta -1).$$
\end{theorem}

The value of $s$ denotes the number of sieving primes $p_1,p_2,p_3,\ldots ,p_s$ which divide $p-1$ but not $e$, an integer of our own choosing. This $e$ is an even divisor of $p-1$ and must be chosen so that $\delta>0$ and $(2^{n-s}\Delta -1)$ is minimised. Once this $e$ is chosen to give our optimal value of $s$, we can proceed to solve the inequality for $p$.

Firstly, for each value of $\omega(p-1)$, $p-1=\prod_{i \leq \omega(p-1)}{p_i}$, the product of the first $\omega(p-1)$ primes, is checked to see if it satisfies the inequality. For values of $p-1=\prod_{i \leq \omega(p-1)}{p_i}<2.5 \times 10^{15}$, the larger value is used as we know all primes below $2.5 \times 10^{15}$ satisfy the values of $\alpha$ being tested. The values of $\omega(p-1)$ which still do not satisfy this equality must be tested further using the prime divisor tree.

For each of these, the value of $p$ which does make the inequality true is set as an upper bound $p_u$, leaving all $p \in[{2.5 \times 10^{15}, p_u}]$ as exceptions which must be checked to see if they have sufficiently small primitive roots. After use of the sieve, we have proven \ref{single}. $\alpha=0.6309$ holds for all values of $\omega(p-1) <22$, and the minimum value of $p$ with $\omega(p-1)=23$ is $p=2.67 \times 10^{32}$. Applying the prime divisor tree (described in \S6) to $\alpha=0.6309$ generated too many exceptions for the computation to be completed on a standard laptop, and so no further investigation into $\alpha=0.6309$ was completed.

\begin{table}[h]
\begin{center}
\caption{Exceptions of $\omega (p-1)$ after use of sieve}

 \begin{tabular}{||c c c||}
 \hline
 $\alpha$ & Lower bound on $w(n)$  & Upper bound on $w(n)$ \\ [0.5ex]
 \hline\hline
 0.69 & - & - \\
 \hline
 0.68 & 13 & 13 \\
 \hline
 0.65 & 5 & 18 \\
 \hline
 0.6309 & 5 & 22 \\ [1ex]
 \hline
\end{tabular}
\end{center}
\end{table}

\section{Computation: the prime divisor tree}

The prime divisor tree relies on the same sieve developed in \S4, but recalculates the optimal $s$ and $\delta$ at each node on a tree. This tree splits up the remaining primes $p$ to be checked according to the specific primes dividing $p-1$. The program was run in SageMath, based on code developed by McGown, Trevi\~no and Trudgian \cite{MTT}, and Hunter \cite{H}.

The prime divisor tree branches off from each point according to whether the next sequential prime divides $p-1$. Since we know $2$ already divides, the first two branches of the tree split up primes according to whether $3$ does or does not divide $p-1$. Each of these nodes may branch further according to whether $5$ does or does not divide $p-1$, and so on.

At the nodes where the prime does not divide $p-1$, this prime can be removed from the calculation of $\delta$, and the next sequential prime, the $\omega(p-1)+1$th prime is included. This increases $\delta$ and so the upper bound on the primes which must still be checked is reduced. In some cases, the upper and lower bounds overlap, and so all possible exceptions are eliminated. In this case, this particular branch of the tree terminates and no further nodes are introduced. If the bounds do not overlap, we check how many prime numbers are in the new reduced interval, and hence how many times we need to check for a primitive root.

At nodes where the prime in question \textit{does} divide $p-1$, we know that all primes must be of the form $(p-1)= k \times m$ where $k$ is the product of those primes that \textit{do} divide $p-1$, and $m$ is the product of the still unknown prime divisors of $p-1$. For example, if we know that $2,3$ and $5$ all divide $p-1$, then $k= 2 \times 3 \times 5 = 30$ Since the first few prime divisors are now known, this eliminates many of the primes within the interval determined by the sieve, reducing the number of primes where we must search for a small primitive root.

Through this process, each branches reduces the number of primes to be checked. If the number to be checked is sufficiently small ($<10^{5}$), the primes are enumerated and primitive roots found using the inbuilt primitive root finder in SageMath. If the interval of exceptions is still too large, the tree branches further until all exceptions have been enumerated.

Through the use of the prime divisor tree, the value of $\alpha$ that holds for all $p$ can be reduced slightly to $0.68$. This improvement is the last step in establishing Theorem 2. Like Hunter, I found that the number of exceptions to be checked increases very rapidly with small reductions in $\alpha$ below $0.68$, and so run times for the code expand to impractical lengths. Below $\alpha=0.68$, my computer unable to run the code to completion, and hence could not lower the value of $\alpha$ further.

\section{Future work}
Through the combined use of the P\'olya--Vinogradov and Burgess bounds, we have improved the bounds on the least primitive root modulo a prime. Specifically, we have reduced the universal bound on the smallest primitive roots for all primes to $g(p) = p^{0.68}$ (Theorem 2). We have also lowered the minimum prime for which $g(p)=p^{0.6309}$ holds, to $p=2.67 \times 10^{32}$ (Theorem 1). While this marks a significant improvement on previous results, as discussed above, it leaves room for future work to lower the bound further.

The largest improvements are likely to come from a tighter constant for the P\'olya--Vinogradov Inequality, as this will increase the strength of the sieve. Due to the number of computations required to run the prime divisor tree, it is only capable of eliminating $2-4$ value of $\omega(p-1)$ in a reasonable time frame. It is therefore preferable to eliminate as many values of $\omega(p-1)$ as possible before we reach the prime divisor tree. One possible source of improvement would be to develop a means of splitting the Dirichlet character sum into odd and even cases. As outlined in Theorem 2 of \cite{FS}, different bounds exist for odd and even Dirichlet characters. In this paper, the weaker bound was taken for all characters as we had no means to systematically split the character sum.

Running the code on a more powerful computer would likely also deliver some benefit. If the values of $\omega(p-1)$ that could be easily checked could be increased, proving that $\alpha = 0.6309$ holds for all $p$ could become feasible.
\section{Acknowledgments}

I'd like to thank Dr Tim Trudgian for his help developing a research topic and providing invaluable guidance as my supervisor. I'd also like to thank UNSW Canberra for providing funding through the Summer Scholarship program and Dr Kevin McGowan for providing feedback on the Sage code.


\begin{thebibliography}{9}
\bibitem{B} D.A. Burgess, \textit{On character sums and primitive roots}, Proc. Lond. Math. Soc. 12(3) (1962), 179-192.
\bibitem{COT} S.D. Cohen, T. Oliveira e Silva, T. Trudgian, \textit{On Grosswald's conjecture on primitive roots}, Acta Arith., 172(3) (2016), 263-270.
\bibitem{CST} S.D. Cohen and T. Trudgian, \textit{On the least square-free primitive root modulo p}, J. Number Theory, 170 (2017), 10-16.
\bibitem{FS} D.A. Frolenkov and K. Soundararajan, \textit{A generalization of the P\'olya–Vinogradov inequality}, Ramanujan J. 31(3) (2013), 271-279.
\bibitem{H} M. Hunter, \textit{The Least Square-free Primitive Root Modulo a Prime}, ANU, Honours Thesis, (2016).
\bibitem{MTT} K. McGown, E. Trevi\~no, and T. Trudgian, \textit{Resolving Grosswald’s conjecture on GRH}, Funct. Approx. Comment. Math. 55(2), (2016), 215-225.
\bibitem{R} G. Robin, \textit{Estimation de la fonction de Tchebychef} $\theta$ \textit{sur le k-i\`eme nombre premier et grandes valeurs de la fonction} $\omega(p-1)$ \textit{nombre de diviseurs premiers de n}, Acta Arith. 42 (1983), 367-369.
\bibitem{T} E. Trevi\~no, \textit{The Burgess inequality and the least} k\textit{th power non-residue}, Intl. Journal of Number Theory 11(5) (2015) 1653-1678.
\end{thebibliography}
\end{document}